\def\ver{June 19, 2001, v.2}
\documentstyle{amsppt}
\magnification=1200
\hsize=6.5truein
\vsize=8.9truein
\topmatter
\title On Mildenhall's Theorem
\endtitle
\author Morihiko Saito
\endauthor
\affil RIMS Kyoto University, Kyoto 606-8502 Japan \endaffil
\abstract We show that Mildenhall's theorem implies that the
indecomposable higher Chow group of a self-product of an elliptic
curve over the complex number field is infinite dimensional, if the
elliptic curve is modular and defined over rational numbers.
For the moment we cannot prove even the nontriviality of the
indecomposable higher Chow group by a complex analytic
method in this case.
\endabstract
\endtopmatter
\tolerance=1000
\baselineskip=12pt
\def\bC{{\Bbb C}}
\def\bF{{\Bbb F}}
\def\bP{{\Bbb P}}
\def\bQ{{\Bbb Q}}
\def\bZ{{\Bbb Z}}
\def\cE{{\Cal E}}
\def\cK{{\Cal K}}
\def\cO{{\Cal O}}
\def\cX{{\Cal X}}
\def\Pic{\hbox{{\rm Pic}}}
\def\div{\hbox{{\rm div}}\,}
\def\NS{\hbox{{\rm NS}}}
\def\CH{\hbox{{\rm CH}}}
\def\pr{\hbox{{\rm pr}}}
\def\Spec{\hbox{{\rm Spec}}\,}
\def\End{\hbox{{\rm End}}}
\def\ind{\text{\rm ind}}
\def\dec{\text{\rm dec}}
\def\gind{\text{\rm g-ind}}
\def\gdec{\text{\rm g-dec}}
\def\GCM{\text{\rm GCM}}
\def\simto{\buildrel\sim\over\to}
\def\SameAuthor{\vrule height3pt depth-2.5pt width1cm}

\document
\centerline{{\bf Introduction}}
\bigskip
\noindent
Let
$ X $ be a smooth algebraic variety.
Consider a finite number of rational functions
$ g_{j} $ on irreducible closed subvarieties
$ Z_{j} $ of codimension
$ p $ in
$ X $ such that
$ \sum \div g_{j} = 0 $ on
$ X $.
This naturally appears, for example, when we study the kernel of
$ \CH^{p}(Z) \to \CH^{p}(Y) $ for a closed subvariety
$ Z $ of
$ Y $ with
$ X = Y \setminus Z $.
Such objects are called higher cycles.
Modulo a suitable equivalence relation defined by tame symbols, they form
a group which is denoted by
$ H^{p}(X,\cK_{p+1}) $.
Here
$ \cK_{p+1} $ is the Zariski-sheafification of the Quillen
$ K $-group [15], and the relation with the above description follows from
the
Gersten resolution.
It is known (see e.g. [13]) that this group is naturally isomorphic to the
higher Chow group
$ \CH^{p+1}(X,1) $,
and the above relation to the Chow group is explained also by the
localization sequence [4].

If the
$ g_{j} $ are constant, a higher cycle
$ \sum (Z_{j},g_{j}) $ is considered rather trivial, and is called
{\it decomposable}.
The subgroup consisting of such cycles is denoted by
$ \CH_{\dec}^{p+1}(X,1) $.
We define the group of {\it indecomposable} higher cycles
$ \CH_{\ind}^{p+1}(X,1) $ to be its quotient group.
Since we often neglect torsion, we usually consider
$ \CH_{\ind}^{p+1}(X,1)_{\bQ} $ where the subscript of
$ \bQ $ means the tensor with
$ \bQ $.
If
$ X $ is a variety over a field
$ k $ which is not algebraically closed, it is sometimes useful to
consider {\it geometrically decomposable} higher cycles, which are
cycles whose base change by a finite extension of
$ k $ is decomposable.
We can define similarly the groups of geometrically decomposable (or
indecomposable) higher cycles
$ \CH_{\gdec}^{p+1}(X,1) $ (or
$ \CH_{\gind}^{p+1}(X,1)) $.

It is often observed that constructing a nonzero element of
$ \CH_{\ind}^{p+1}(X,1)_{\bQ} $ for a given variety
$ X $ is not an easy task.
Most of the constructions over the complex number field consider a family
of varieties with higher cycles, and show that for a general member the
cycle is indecomposable (see e.g. [6], [7], [8], [10], [13], [17], [21]).
However, the situation is different over a number field.
Mildenhall proved

\medskip\noindent
{\bf 0.1.~Theorem} ([14]). {\it
Let
$ E $ be a modular elliptic curve over
$ \bQ $, and
$ X $ the self-product of
$ E $.
Then
$ \CH_{\ind}^{2}(X,1)_{\bQ} $ is infinite dimensional.
If furthermore
$ E $ has complex multiplication such that
$ \End_{K}(E)_{\bQ} $ is a quadratic imaginary field
$ K $, then
$ \CH_{\ind}^{2}(X_{K},1)_{\bQ} $ is also infinite dimensional.
}

\medskip
If
$ E $ has complex multiplication, let
$ K $ be as above, and
$ K = \bQ $ otherwise.
Let
$ S = \Spec \bZ[1/6N] $ where
$ N $ is the discriminant of
$ E $.
Then he showed

\medskip\noindent
{\bf 0.2.~Theorem} ([14]). {\it
The divisor map
$ \div : \CH^{2}(X,1)_{\bQ} \to \oplus_{s\in |S|}
\End_{k(s)}(E_{s})_{\bQ} $ is surjective.
}

\medskip
Here
$ |S| $ is the set of the closed points of
$ S $,
$ \End_{k(s)}(E_{s}) $ is identified with a quotient of the N\'eron-Severi
group
$ \NS(X_{s}) $,
and the divisor map is defined by extending a higher cycle to the abelian
scheme over
$ S $ whose generic fiber is
$ X $.
The proof uses modular units and
$ p $-Hecke correspondence to show that a multiple of the
$ p $-th power Frobenius belongs to the image of the divisor map,
see also [9], [18].
Unfortunately, it is not explained in [14] that Theorem (0.2) implies
more than (0.1).
In fact, we can deduce from (0.2)

\medskip\noindent
{\bf 0.3.~Theorem.} {\it
Let
$ E, X $ and
$ K $ be as above.
Then
$ \CH_{\gind}^{2}(X_{K},1)_{\bQ} $ and
$ \CH_{\ind}^{2}(X_{\bC},1)_{\bQ} $ are infinite dimensional.
}

\medskip
We can omit the hypothesis on the modularity of
$ E $ in most cases due to Wiles [22] and others.
In the CM case, this is classically well-known (due to Deuring).
The proof of (0.3) uses a model of
$ X $,
and is rather standard (see [2], [3], [16], [20]).
It should be noted that we cannot prove the nonvanishing of
$ \CH_{\ind}^{2}(X_{\bC},1)_{\bQ} $ by a complex analytic method for
the moment.
For example, the image of the cycle of Gordon and Lewis [10] in Deligne
cohomology vanishes if
$ E $ has complex multiplication.

I would like to thank M\"uller-Stach for drawing my attention to
Mildenhall's paper [14] and also for useful discussions.

\bigskip\bigskip\centerline{{\bf 1. Proof of Theorem (0.3)}}

\bigskip
\noindent
{\bf 1.1.}
Let
$ E, X, K, N $ and
$ S $ be as in the introduction.
Let
$ \cE $ be the abelian scheme on
$ S $ whose generic fiber is
$ E $.
Let
$ \cX = \cE\times_{S}\cE $ so that the generic fiber of
$ \cX $ is
$ X $.
Recall that
$ K $ is a quadratic imaginary field if
$ E $ is CM-type, and
$ K = \bQ $ otherwise.
For a number field
$ L $ containing
$ K $, let
$ S_{L} = \Spec \cO_{L}[1/6N] $, where
$ \cO_{L} $ is the ring of integers of
$ L $.
We have
$$
\End_{L}(E_{L})_{\bQ} = K,
\leqno(1.1.1)
$$
(see e.g. [19], II.2.2), where
$ E_{L} = E\otimes_{\bQ}L $.
Note that
$$
\End_{\bQ}(E) = \bZ.
\leqno(1.1.2)
$$
Let
$ \cE_{L} = \cE\otimes_{S}S_{L} $.
It is an abelian scheme over
$ S_{L} $, and is the N\'eron model of
$ E_{L} $ (see [5]).
So the restriction induces an isomorphism
$$
\End_{S_{L}}(\cE_{L}) = \End_{L}(E_{L}).
\leqno(1.1.3)
$$

For
$ s \in |S_{L}| $,
let
$ k(s) $ denote the residue field of
$ s $,
and
$ E_{s}, X_{s} $ the fiber of
$ \cE_{L}, \cX_{L} $ at
$ s $.
If
$ s' \in |S| $ underlies
$ s $, then
$ E_{s} $ is the base change of
$ E_{s'} $ by
$ k(s') \to k(s) $ (and similarly for
$ X_{s} $).
For
$ s \in |S_{L}| $, the restriction map together with (1.1.3) induces
$$
\End_{L}(E_{L})_{\bQ} \to \End_{k(s)}(E_{s})_{\bQ}.
\leqno(1.1.4)
$$
We denote its image by
$ \End_{k(s)}(E_{s})^{\GCM}_{\bQ} $ (GCM for global complex
multiplication).
We have a canonical isomorphism
$$
\End_{k(s)}(E_{s}) = \Pic(E_{s}\times_{k(s)}E_{s})/(\pr_{1}^{*}
\Pic(E_{s}) + \pr_{2}^{*}\Pic(E_{s}))
\leqno(1.1.5)
$$
by the theory of algebraic correspondences.
In particular,
$ \End_{k(s)}(E_{s}) $ is identified with a quotient group of the
N\'eron-Severi group
$ \NS(E_{s}) $.

\medskip\noindent
{\bf 1.2.~Lemma.} {\it
For
$ s \in |S_{L}| $ over
$ s' \in |S_{K}| $,
the push-forward of cycles induces}
$$
\End_{k(s)}(E_{s})_{\bQ}^{\GCM} \to
\End_{k(s')}(E_{s'})_{\bQ}^{\GCM}.
\leqno(1.2.1)
$$

\medskip\noindent
{\it Proof.} The pull-back of cycles induces an isomorphism
$$
\End_{k(s')}(E_{s'})_{\bQ}^{\GCM} \simto
\End_{k(s)}(E_{s})_{\bQ}^{\GCM},
\leqno(1.2.2)
$$
using the commutative diagram
$$
\CD
\End_{S_{L}}(\cE_{S_{L}}) @>>> \End_{k(s)}(E_{s})
\\
@AAA @AAA
\\
\End_{S_{K}}(\cE_{S_{K}}) @>>> \End_{k(s')}(E_{s'}),
\endCD
\leqno(1.2.3)
$$
because the left vertical morphism is an isomorphism.
Then the assertion follows from the fact that the composition of
the pull-back and the push-down of cycles under a finite \'etale morphism
is the multiplication of the identity by the degree.

\medskip\noindent
{\bf 1.3.}
{\bf Proof of (0.3).} Since the composition of the pull-back and the
push-down of higher cycles under the morphism
$ X_{L} \to X_{K} $ is the identity multiplied by
$ [L:K] $, and the divisor map is compatible with the push-down of cycles,
we see that the images of geometrically decomposable cycles by the divisor
map
$$
\div : \CH^{2}(X_{K},1)_{\bQ} \to \oplus_{s\in |S_{K}|}
\End_{k(s)}(E_{s})_{\bQ}
\leqno(1.3.1)
$$
are contained in
$ \oplus_{s\in |S_{K}|} \End_{k(s)}(E_{s})_{\bQ}^{\GCM} $
due to (1.2).
In the CM case, it is known that there are infinitely many
$ s \in |S_{K}| $ such that
$ k(s) = \bF_{p^{2}} $ and
$ E_{s} $ is a supersingular elliptic curve, i.e.
$ \dim_{\bQ} \End_{k(s)}(E_{s})_{\bQ} = 4 $, see e.g. [11], [14].
Indeed, if
$ (D/p) = -1 $ with
$ K = \bQ(\sqrt{D}) $, then a prime number
$ p $,  which is identified with
$ s' \in |S| $, does not split in
$ K $,  and
$ E_{s'} $ is supersingular, see e.g. [11], 13, Th. 5.
The existence of infinitely many such
$ p $ follows from the Tchebotarev density theorem or the Dirichlet
theorem together with the quadratic reciprocity law.
In this case, the characteristic polynomial of the
$ p $-th power (geometric) Frobenius is
$ T^{2} + p $ as well-known, so that the action of the
$ p^{2} $-th power Frobenius is
$ - p $, and the
$ p $-th power Frobenius is not contained in
$ \End_{k(s)}(E_{s})_{\bQ}^{\GCM} $.
In general,
$ \dim_{\bQ} \End_{k(s)}(E_{s}) \ge 2 $ by the Frobenius.
So Theorem (0.2) implies
$$
\dim_{\bQ} \CH_{\gind}^{2}(X_{K},1)_{\bQ} = \infty.
\leqno(1.3.2)
$$

From this we can deduce
$$
\dim_{\bQ} \CH_{\ind}^{2}(X_{\bC},1)_{\bQ} = \infty ,
\leqno(1.3.3)
$$
if we know the injectivity of
$$
\CH_{\gind}^{2}(X_{K},1) \to \CH_{\ind}^{2}(X_{\bC},1).
\leqno(1.3.4)
$$
To show this injectivity, let
$ \zeta \in \CH^{2}(X_{K},1) $,
and assume
$ \zeta \otimes_{K}\bC \in \CH^{2}(X_{\bC},1) $ is decomposable.
Then there exist a finitely generated
$ K $-subalgebra
$ R $ of
$ \bC $,
divisors
$ Z_{i} $ on
$ X\otimes_{K}R, \alpha_{i} \in R \,(1 \le i \le a) $ and rational
functions
$ f_{j}, g_{j} $ on
$ X\otimes_{K}R \,(1 \le i \le b) $ such that
$$
\pr_{1}^{*}\zeta - \sum_{i} (Z_{i},\alpha_{i}) = \sum_{j} \{f_{j},g_{j}\},
\leqno(1.3.4)
$$
where
$ \{f_{j},g_{j}\} $ denotes the tame symbol.
Restricting to a general closed point of
$ \Spec R $,
we may assume
$ R $ is a finite extension of
$ K $ by replacing
$ Z_{i}, \alpha_{i}, f_{j}, g_{j} $ with their restrictions to the fiber
over the closed point.
Then (1.3.4) means that
$ \zeta $ is geometrically decomposable, and the assertion follows.

\medskip\noindent
{\bf 1.4.}
{\bf Remark.} It does not seem easy to construct explicitly an
indecomposable higher cycle on a self-product of an elliptic curve of CM
type.
For example, we can show the vanishing of the image of the cycle of
Gordon-Lewis [10] by the transcendental part of the Abel-Jacobi map (i.e.
by the integral along a non closed
$ C^{\infty } $ chain, see [12]).
Indeed, the cycle depends actually on the choice of the function
$ f $ which comes from a rational function on
$ \bP^{1} $ (underlying
$ E_{1}), $ and we can move the zero of this function toward the pole so
that the above
$ C^{\infty } $ chain becomes smaller and smaller, but the image stays
constant due to the rigidity of Beilinson [1] and M\"uller-Stach [13].
However it seems also difficult to show that this higher cycle is really
decomposable.
This is related to the problem on the injectivity of the reduced
Abel-Jacobi map (see e.g. [17]).

\newpage\centerline{{\bf References}}\bigskip

\item{[1]}
A. Beilinson, Higher regulators and values of
$ L $-functions, J. Soviet Math. 30 (1985), 2036--2070.

\item{[2]}
A. Beilinson, J. Bernstein and P. Deligne, Faisceaux pervers,
Ast\'erisque, vol. 100, Soc. Math. France, Paris, 1982.

\item{[3]}
S. Bloch., Lectures on algebraic cycles, Duke University Mathematical
series 4, Durham, 1980.

\item{[4]}
\SameAuthor, Algebraic cycles and higher K-theory, Advances in Math., 61
(1986), 267--304.

\item{[5]}
S. Bosch, W. L\"utkebohmert and M. Raynaud, N\'eron models, Springer,
Berlin, 1990.

\item{[6]}
A. Collino, Griffiths' infinitesimal invariant and higher
$ K $-theory on hyperelliptic Jacobians, J. Alg. Geom. 6 (1997), 393--415.

\item{[7]}
A. Collino and N. Fakhruddin, Indecomposable higher Chow cycles on
Jacobians, preprint.

\item{[8]}
P. del Angel and S. M\"uller-Stach, The transcendental part of the
regulator map for
$ K_{1} $ on a mirror family of
$ K3 $ surfaces, preprint.

\item{[9]}
M. Flach, A finite theorem for the symmetric square of an elliptic curve,
Inv. Math. 109 (1992), 307--327.

\item{[10]}
B. Gordon and J. Lewis, Indecomposable higher Chow cycles on products
of elliptic curves, J. Alg. Geom. 8 (1999), 543--567.

\item{[11]}
S. Lang, Elliptic functions, Springer, Berlin, 1987.

\item{[12]}
M. Levine, Localization on singular varieties, Inv. Math. 91 (1988),
423--464.

\item{[13]}
S. M\"uller-Stach, Constructing indecomposable motivic cohomology
classes on algebraic surfaces, J. Alg. Geom. 6 (1997), 513--543.

\item{[14]}
S.J.M. Mildenhall, Cycles in a product of elliptic curves, and a group
analogous to the class group, Duke Math. J. 67 (1992), 387--406.

\item{[15]}
D. Quillen, Higher algebraic
$ K $-theory I, Lect. Notes in Math., vol. 341, Springer, Berlin, 1973,
pp. 85--147.

\item{[16]}
M. Saito, Arithmetic mixed sheaves, preprint (math.AG/9907189).

\item{[17]}
\SameAuthor, Bloch's conjecture, Deligne cohomology and higher Chow
groups, preprint RIMS--1284.

\item{[18]}
A.J. Scholl, On Modular units, Math. Ann. 285 (1989), 503--510.

\item{[19]}
J.H. Silverman, Advanced topics in the arithmetic of elliptic curves,
Springer, Berlin, 1994.

\item{[20]}
C. Voisin, Transcendental methods in the study of algebraic cycles, in
Lect. Notes in Math. vol. 1594, pp. 153--222.

\item{[21]}
\SameAuthor, Nori's connectivity theorem and higher Chow groups, preprint.

\item{[22]}
A. Wiles, Modular elliptic curves and Fermat's last theorem, Ann. Math.
142 (1995), 443--551.

\bigskip
\noindent
\ver

\bye